\documentclass[number,citesort,dvips]{arxbj}

\aid{0}
\volume{18}
\issue{1}
\pubyear{2012}
\firstpage{1}
\lastpage{23}
\doi{10.3150/10-BEJ324}

\makeatletter
\newtheorem{theo}{Theorem}
\newremark{remi}{Remark}
\newtheorem{lemme}[theo]{Lemma}
\newtheorem{prop}[theo]{Proposition}

\newcommand{\eqref}[1]{(\ref{#1})}

\newcommand{\1}{\mathbf{1}}
\newcommand{\R}{\mathbb R}
\newcommand{\E}{\mathbb E}
\newcommand{\oT}{[0{,}T]}
\newcommand{\eps}{\varepsilon}
\newcommand{\rd}{\mathbb R^{d}}
\newcommand{\s}{\sigma}
\newcommand{\ffi}{\varphi}

\makeatother

\begin{document}
\begin{frontmatter}

\title{Transportation inequalities for stochastic differential
equations driven by a fractional Brownian motion}
\runtitle{Transportation inequalities for SDEs driven by a fBm}

\begin{aug}
\author{\fnms{Bruno} \snm{Saussereau}\corref{}\ead[label=e1]{bruno.saussereau@univ-fcomte.fr}}
\runauthor{B. Saussereau}
\address{Laboratoire de Math\'{e}matiques de Besan\c{c}on, CNRS, UMR 6623,
16 Route de Gray,
25030 Besan\c{c}on cedex, France. \printead{e1}}
\end{aug}

% HISTORY:
\received{\smonth{2} \syear{2008}}
\revised{\smonth{8} \syear{2010}}

% ABSTRACT
%
\begin{abstract}
We establish Talagrand's $T_1$ and $T_2$ inequalities for the law of
the solution of a stochastic differential equation driven by a
fractional Brownian motion with Hurst parameter $H>1/2$. We use the
$L^2$ metric and the uniform metric on the path space of continuous
functions on $[0,T]$. These results are applied to study small-time and
large-time asymptotics for the solutions of such equations by means of
a Hoeffding-type inequality.
\end{abstract}

% KEYWORDS
%

\end{frontmatter}
%

%s1 ###
\section{Introduction}
Suppose that $B^H=(B_t^H)_{t\in[0,T]}$ is an $m$-dimensional fractional
Brownian motion (fBm) with Hurst parameter $H$ defined on a complete
filtered probability space $(\Omega, \mathcal F, (\mathcal
F_t)_{t\in[0,T]},\mathbb P)$. By this, we mean that the components
$B^{H,j}$, $j=1,\ldots,m$, are independent centered Gaussian processes
with the covariance function
\[
R_H(s,t)=\tfrac{1}2( t^{2H}+s^{2H}-|t-s|^{2H}).
\]
If $H=1/2$, then $B^{H}$ is clearly a Brownian motion. Since for any
$p\ge1$, $ \mathbb E|B^{H,j}_t-B^{H,j}_s|^p=c_p|t-s|^{pH}$, the
processes $B^{H,j}$ have $\alpha$-H\"{o}lder continuous paths for all
$\alpha\in(0,H)$ (see~\cite{n} for further information about fBm).

In this article we fix $1/2< H<1$ and are interested in the solution
$(X_t)_{t\in[0,T]}$ of the stochastic differential equation
\begin{equation}\label{eds}
 X^{i}_t = x^{i}+\sum_{j=1}^{m}\int_{0}^t \s^{i,j}(X_s)\,\mathrm{d}
B^{H,j}_s + \int_{0}^t b^i(X_s)\,\mathrm{d}s,\qquad t\in\oT ,
\end{equation}
$i=1,\ldots,d$, where $x\in\rd$ is the initial value of the
process\vadjust{\goodbreak}
$X$.

Under suitable assumptions on $\s$, the processes % $(\s(X_s))_{s\in
$\sigma(X)$ and $B^H$ have trajectories which are H\"{o}lder continuous
of order strictly larger than $1/2,$ so we can use the integral
introduced by Young in \cite{young}. The stochastic integral in
\eqref{eds} is then a pathwise Riemann--Stieltjes integral. A first
result on the existence and uniqueness of a solution of such an
equation was obtained in \cite{lyons} using the notion of
$p$-variation. The theory of rough paths introduced by Lyons in
\cite{lyons} was used by Coutin and Qian in order to prove an existence
and uniqueness result for the equation \ref{eds} (see~\cite{cq}). The
Riemann--Stieltjes integral appearing in equation (\ref{eds}) can be
expressed as a Lebesgue integral using a fractional integration by
parts formula (see Z\"{a}hle \cite{za}). Using this formula, Nualart
and R\u{a}\c{s}canu have established in \cite{nr} the existence of a
unique solution for a class of general differential equations that
includes (\ref{eds}). Regularity (in the sense of Malliavin calculus)
and absolute continuity of the law of the random variables $X_t$ have
since been investigated in \cite{bh,hn,nsimon,ns}.

This work is strongly motivated by the study of the small-time and
large-time behaviors of the solution of \eqref{eds}. To the best of our
knowledge, little seems to be known on this subject. In~\cite{h} the
author investigates the ergodicity of the solution when $\sigma$ is
constant, as well as the convergence rate toward the stationary
solution; see also \cite{ms} for infinite-dimensional evolution
equations driven by an fBm in an additive way. We will be able to state
small-time and large-time asymptotic properties as consequences of
stronger properties: the concentration inequalities on the path space
of continuous functions.

For several years, the transportation cost-information inequalities and
their applications to diffusion processes have been widely studied. In
this paper we apply recent results on fractional differential equations
in order to obtain Talagrand's inequalities. Let us now consider the
kinds of inequalities we will deal with. To measure distances between
probability measures, we use transportation distances, also called
Wasserstein distances. Let $(E,d)$ be a metric space equipped with a
$\sigma$-field $\mathcal B$ such that the distance~$d$ is $\mathcal
B\otimes\mathcal B$-measurable. Given $p\in[1,+\infty]$ and two
probability measures $\mu$ and $\nu$ on $E$, the Wasserstein distance
is defined by\vspace*{-2pt}
\[
W_p^d(\mu,\nu) = \inf\biggl( \int\!\!\!\int
\mathrm{d}(x,y)^p\,\mathrm{d}\pi(x,y) \biggr)^{1/p},\vspace*{-2pt}
\]
where the infimum is taken over all the probability measures $\pi$ on
$E\times E$ with marginal distributions $\mu$ and $\nu$. The relative
entropy of $\nu$ with respect to $\mu$ is defined as\vspace*{-2pt}
\[
\mathbf{H}(\nu/\mu) = \cases{
\displaystyle\int\log\frac{\mathrm{d}\nu}{\mathrm{d}\mu}\,\mathrm{d}\nu, &\quad if $\nu\ll\mu$, \cr
+\infty, &\quad otherwise.}\vspace*{-2pt}
\]
The probability measure $\mu$ satisfies the $L^p$ transportation
inequality on $(E,d)$ if there exists a~constant $C\ge0$ such that for
any probability measure $\nu$,\vspace*{-2pt}
\[
W_p^d(\mu,\nu) \le\sqrt{2C \mathbf{H}(\nu/\mu)}.\vspace*{-2pt}
\]
As usual, we write $\mu\in T_p(C)$ for this relation. The properties
$T_1(C)$ and $T_2(C)$ are of particular interest. The phenomenon of
measure concentration is related to $T_1(C)$ (see the monograph of
Ledoux \cite{ledoux}).

The property $T_2(C)$ is stronger than $T_1(C)$ but is not so well
characterized. It was first established by\vadjust{\goodbreak} Talagrand \cite{tala} for
the Gaussian measure and generalized in \cite{fu2} to the framework of
an abstract Wiener space; see \cite{bgl,ov} for the relationship
between $T_2(C)$ and other properties such as the Poincar\'{e}
inequality and Hamilton--Jacobi equations. The logarithmic Sobolev
inequality introduced by Gross \cite{gross} plays a particular role in
this theory since it implies $T_2(C)$ (see \cite{bgl,ov,wang2}).

With regard to the paths of diffusion processes, the $T_2$
transportation inequality with respect to the Cameron--Martin metric
was proven in \cite{dgw} by means of the Girsanov transform. The
authors also provided a direct proof of the $T_1$ trans\-por\-tation
inequality with respect to the uniform metric using the Gaussian tail
criterion (see Section \ref{surprise} for more details). Later, in
\cite{wz}, $T_2(C)$ was established with respect to the uniform
metric. Finally, Gourcy and Wu \cite{gw} established the log-Sobolev
inequality for the Brownian motion with drift in the $L^2$ metric
instead of the usual Cameron--Martin metric. As a consequence, they
derived the $T_2(C)$ property with respect to this metric and a
concentration inequality (of correct order for large time) for some
functionals of the process. In \cite{wang1}, the~$T_2(C)$ property with
respect to the $L^2$ metric was established for elliptic diffusions on
a~Riemannian manifold.

In this paper, we investigate the properties $T_1(C)$ and $T_2(C)$ for
the law $\mathbb P_x$ of the solution $(X_t)_{0\le t\le T}$ of the
equation \eqref{eds} in various situations. We work on the space of
continuous functions endowed with the uniform metric or the $L^2$
metric. $T_2(C)$ will hold for a multidimensional equation when $\sigma
= I_d$ and $d=m,$ and for a one-dimensional equation when the diffusion
coefficient $\sigma$ is non-constant. It will also be established with
respect to the uniform distance rather than the $L^2$ metric. The use
of this second metric will be of particular interest when dealing with
large-time asymptotics. The~$T_1(C)$ property will be proven for a
multidimensional equation with a diffusion matrix $\sigma$ that is only
a time-dependent function. In the one-dimensional case, the function~$\sigma$
may depend on the space variable. This property is proved with
respect to the uniform metric for small-time horizon $T$. This
restriction to small time is discussed after Theorem \ref{t1-dim1} and
this result is of great interest when we apply it to small-time
asymptotics.

The paper is organized as follows. Section \ref{nota} is devoted to the
statement of our results. In Section \ref{appli}, we review the usual
consequences of transportation inequalities for large- and small-time
behavior. Section \ref{prelim} contains the estimation of the
difference of the solutions of two deterministic differential equations
driven by H\"{o}lder continuous functions of order greater than $1/2$.
The method we develop to prove our main results in Section
\ref{preuves} is the counterpart of the usual case: the Gaussian
integrability condition for $T_1(C)$, Girsanov's formula and an
explicit control for a~specific coupling of two paths of the solution
of the stochastic differential equation. In the framework of fractional
Brownian motion, this control is new, to the best of our knowledge. In
Section~\ref{surprise} we make a quite surprising remark about the link
between the constant $C$ in a property $T_1(C)$ and a Gaussian tail. A
priori this remark is independent of the rest of this work, but it can
be helpful when trying to prove $T_1(C)$ via an exponential moment.
Finally, a Fernique-type lemma is proved in the \hyperref[appendix]{Appendix}.

%
%The paper is organized as follows. First we give the framework we
%use and state our results. Then we recall and establish some results
%about
%regularity of differential equations driven by H\"{o}lder continuous
%functions of order greater than $1/2$.
%
%Finally, some technical lemmas are eventually proved in an appendix.
%
%s2 ###
\vspace*{-3pt}\section{Main results}\vspace*{-3pt}\label{nota}
We consider a complete probability space $(\Omega, \mathcal F, \mathbb
P)$ on which an $m$-dimensional Brownian motion $(W_t)_{t\in[0,T]}$
is
defined. We denote by $\mathcal F_t =\sigma(W_s,s\le t)$ the
$\sigma$-field generated by $W$ and\vadjust{\goodbreak} completed with respect to $\mathbb
P$. Finally, $B^H=(B_t^H)_{t\in[0,T]}$ is the $m$-dimensional fBm
defined on $(\Omega,\mathcal F,\mathbb P)$ transferred from $W$. This
means that $B^H$ can be expressed as
%e1 ###
%
\begin{equation}\label{transfer}
B^{H,i}_t =\int_{0}^t K_H(t,s)\,\mathrm{d}W^i_s,\qquad  i=1,\ldots,m,
\end{equation}
where the square-integrable deterministic kernel $K_H$ is defined by
\begin{equation}\label{defK}
K_{H}(t,s)=c_Hs^{ 1/2-H}\int_{s}^t (u-s)^{H-3/2}u^{H- 1/2}\,\mathrm{d}u
\end{equation}
with $ c_H= (\frac{H(2H-1)}{\beta(2-2H,H-1/2 )} )^{ 1/2}$ for
$s<t$ ($\beta$ denotes the beta function). We set $K_H(t,s)=0$ if
$s\geq t$. The process $B^H$ is $\mathcal F_t$-adapted.

We will also need some notation. For $0<\lambda\leq1$ and $0\le a <
b\le T$, we denote by~$C^{\lambda}(a,\allowbreak b;\mathbb{R}^{d})$ the
space of $\lambda$-H\"{o}lder continuous functions $f \dvtx{[a{,}b]}%
\rightarrow\mathbb{R}^{d}$, equipped with the norm
\[
\Vert f\Vert_{\lambda}:=\Vert f\Vert_{a,b,\infty}+\| f\|_{a,b,\lambda},
\]
where
\[
\| f\|_{a,b,\infty} = \sup_{a\leq r\leq b}|f (r)|\quad \mbox{and}\quad
\|f\|_{a,b,\lambda} = \sup_{a\leq r \leq s \leq b}\frac{|f (s)-f
(r)|}{|s-r|^\lambda}.
\]
We simply write $C^\lambda(a,b)$ when $d=1$.

We consider various forms of the stochastic differential equation
\eqref{eds}. We begin with the equation on $\mathbb R^d$
\begin{equation}\label{eqt1-1}
X^{i}_t  = x^{i} + \int_{0}^t b^i(X_s)\,\mathrm{d}s+\sum_{j=1}^{m}\int_{0}^t
\s^{i,j}(s)\,\mathrm{d} B^{H,j}_s,\qquad t\in\oT,\  i=1,\ldots,d
\end{equation}
and make the following assumptions on the coefficients:
\begin{enumerate}
\item[H1(a)] there exists some $L_b$ such that for any $i=1,\ldots,d$ and
any $z,z'\in\mathbb R^d$,
\[
|b(z)-b(z')| \le L_b |z-z'|;
\]
\item[H1(b)] there exists some $\beta>1-H$ such that $\sigma\in
C^\beta(0,T;\mathbb R^{d\times m})$.
% and $\|\sigma''\|_\infty\le c_2$.
% \item
% \item[(H2')] There exists $c_0,c_1$ such
% that $|\s(x)|\leq c_0 + c_1|x|$ and $\|\s'\|_\infty\leq c_1$.
\end{enumerate}
It has been proven in \cite{nr} that under the above assumptions, there
exists a unique adapted stochastic process solution to equation
\eqref{eds} whose trajectories are H\"{o}lder continuous of order
$H-\epsilon$ for any $\epsilon> 0$.

For this kind of equation, we have the following result.

\begin{theo}\label{t1}
Assume that the assumptions $\mathrm{(H1)}$ are satisfied.
Then, for each $0<T\le(2L_b)^{-1}\wedge1$, there exists a
universal
constant $K,$ independent of the initial point $x,$ such that the law
$\mathbb P_x$ of the solution of equation \eqref{eqt1-1} satisfies the
property $T_1(K \|\s\|_{\beta}T^{2H})$ on~$C(0,\allowbreak T;\mathbb{R}^{d}),$ the
space of $\mathbb R^d$-valued continuous functions on $[0,T]$ equipped\vadjust{\goodbreak}
with the metric~$d_{\infty}$ defined by\vspace*{-1pt}
\[
d_\infty(\gamma_1,\gamma_2) = \sup_{0\leq t\leq
T}|\gamma_1(t)-\gamma_2(t)|.\vspace*{-1pt}
%& + \sup_{0\leq s \leq t \leq T}\frac{| ( \gamma_1 (t)-\gamma_1 (s)
% ) - ( \gamma_2 (t)-\gamma_2 (s) ) |}{|t-s|^\beta}  .
\]
\end{theo}
Of course, this result will be useful for small-time asymptotics of the
process $X$. In the one-dimensional case, we will be able, via a
Lamperti transform, to deduce a result for the nonlinear equation\vspace*{-1pt}
\begin{equation}\label{eqdim1}
X_t = x+ \int_0^t b(X_s)\,\mathrm{d}s + \int_0^t \sigma(X_s)\,\mathrm{d} B^H_s  ,\vspace*{-1pt}
\end{equation}{
where the coefficients satisfy:
\begin{enumerate}[H2(a)]
\item[H2(a)] the function $b$ is bounded by $B:=\sup_{x\in\mathbb
R}|b(x)|$ and there exists some $L_b$ such that for any $z,z'\in
\mathbb R$,\vspace*{-1pt}
\[
|b(z)-b(z')| \le L_b|z-z'|;\vspace*{-1pt}
\]
\item[H2(b)] there exist some $\sigma_2>\sigma_1>0$ such that for
any $x\in\mathbb R$,\vspace*{-1pt}
\[
\sigma_1 \le\sigma(x) \le\sigma_2;\vspace*{-1pt}
\]
\item[H2(c)] there exists a constant $L_\sigma$ such that for any
$z,z'\in\mathbb R$,\vspace*{-1pt}
\[
|\sigma(z)-\sigma(z')| \le L_\sigma|z-z'|.\vspace*{-1pt}
\]
\end{enumerate}
\begin{theo}\label{t1-dim1}
Assume that the hypotheses ($\mathrm{H2}$) are satisfied. There
exists a universal constant $K,$ independent of the initial point $x,$
such that the law $\mathbb P_x$ of the solution of equation
\eqref{eqdim1} satisfies the property $T_1(K \sigma_2^2 T^{2H})$ on $
C(0,T;\mathbb{R}),$ provided that $T\le1 \wedge\frac{\sigma_1^2 }{ 2
\sigma_2 ( L_b \sigma_2 + L_\sigma B )}$.
\end{theo}

Before stating the $T_2$ inequalities, we will explain why the
restriction to small time in the statements of the above theorems is in
fact quite natural. Imagine the case where $b=0$ and $d=m=1$. The
processes $X$ and $B^H$ are then equals. It is known (see (\cite{dgw}, Theorem
2.3) or Section \ref{surprise}) that $T_1(C)$ is then equivalent
to the fact that there exists some $\delta>0$ such that\vspace*{-1pt}
\[
C(\delta) = \mathbb E ( \exp\{ \delta \| B^H-\tilde
B^H\|_{0,T,\infty}^2 \}) <\infty,\vspace*{-1pt}
\]
where $B^H$ and $\tilde B^H$ are two independent fractional Brownian
motions. For $f,\tilde f \in C^\beta(0,T)$ with $f(0)=\tilde f(0)$, we
have $\|f-\tilde f\|_{0,T,\infty} \le T^\beta\| f-\tilde
f\|_{0,T,\beta} $. Then,\vspace*{-1pt}
\[
C(\delta) \le\mathbb E ( \exp\{ \delta T^{2\beta} \| B^H-\tilde
B^H\|_{0,T,\beta}^2 \})\vspace*{-1pt}
\]
and with \eqref{esti-exp} from Lemma \ref{lemme-grr} in the \hyperref[appendix]{Appendix},
the above exponential moment will be finite as soon as $\delta
T^{2\beta} \times128 (2T)^{2(H-\beta)} \le1,$ which implies that
$T$
must be small.

We now return to the statements concerning $T_2$ transportation
inequalities. We consider the solution of the stochastic differential
equation \eqref{eqt1-1} and make the following additional stability
assumption on the coefficient $b$:\vadjust{\goodbreak}
\begin{enumerate}
\item[(H3)] There exists some $B\in\mathbb R$ such that for any
$x,y\in\mathbb R^d$,
\[
 \langle x-y,  b(x)-b(y)  \rangle_{\mathbb R^d} \le B |x-y|^2 .
\]
\end{enumerate}
\begin{theo}\label{t2}
We consider $\mathbb P_x$, the law of the solution of the stochastic
differential equation \eqref{eqt1-1}. We assume that
$\mathrm{(H1)}$ and $\mathrm{(H3)}$ are fulfilled.
The probability measure $\mathbb P_x$ satisfies~$T_2(C)$ on the metric
space $C(0,T;\mathbb R^d)$ with:
\begin{enumerate}[(b)]
\item[(a)] $C=(2/|B|)H T^{2H-1}(1\vee e^{(2B+|B|)\times T})\|\sigma\|
_{0,T,\infty}^2 $ with the metric $d_\infty$;
\item[(b)] $C=(2/B^2)HT^{2H-1}\|\sigma\|_{0,T,\infty}^2 c_{B,T}$ with
\[
c_{B,T} := \cases{
\displaystyle\frac{\mathrm{e}^{3BT}-1}{3}, &\quad if $B> 0$,\vspace*{2pt} \cr
1-\mathrm{e}^{BT}, &\quad if $B<0$,
}
\]
when using the metric
\[
d_2(\gamma_1,\gamma_2) = \biggl( \int_0^T |\gamma_1(t)-\gamma_2(t)|^2\,\mathrm{d}t
\biggr)^{1/2}.
\]
\end{enumerate}
\end{theo}

A result for one-dimensional equations with non-constant diffusion
coefficients can be deduced from Theorem \ref{t2}. We assume that
$d=m=1$ and consider the solution of the stochastic differential
equation \eqref{eqdim1}. We make the following assumptions on the
coefficients:
\begin{enumerate}
\item[H4(a)] there exists some $L_b$ such that for any $z,z'\in
\mathbb R$,
\[
|b(z)-b(z')| \le L_b|z-z'|;
\]
\item[H4(b)] there exist some $\sigma_2>\sigma_1>0$ such that for
any $x'\in\mathbb R$,
\[
\sigma_1 \le\sigma(x) \le\sigma_2;
\]
\item[H4(c)] $b$ and $\sigma$ are differentiable,
and there exists some $B\in\mathbb R$ such that for any $x\in\mathbb
R$,
\[
b'(x)\sigma(x) - \sigma'(x)b(x) \le B.
\]
\end{enumerate}
\begin{theo}\label{t3}
Let $d=m=1$ and assume that the assumptions $\mathrm{(H4)}$
hold. The law $\mathbb P_x$ of the solution of the stochastic
differential equation \eqref{eqdim1} then satisfies the property
$T_2(C)$ on the metric space $C(0,T;\mathbb{R})$ where:
\begin{enumerate}[(b)]
\item[(a)] $C=(2\sigma_1\sigma_2^2/|B|)H T^{2H-1}(1\vee
e^{(2B+|B|)\times T/\sigma_1})$ with the metric $d_\infty$;
\item[(b)] $C=(2\sigma_1^2\sigma_2^2/B^2)HT^{2H-1} c_{B,T}$ with
\[
c_{B,T} := \cases{
\displaystyle\frac{\mathrm{e}^{3BT/\sigma_1}-1}{3}, &\quad if $B> 0$,\vspace*{2pt} \cr
1-\mathrm{e}^{BT/\sigma_1}, &\quad if $B<0$,
}
\]
when one uses the metric $d_2$.\vadjust{\goodbreak}
\end{enumerate}
\end{theo}

We note that (H4) implies (H3) when $d=m=1$ and
$\sigma$ is identically equal to $1$.

The constants $C$ in the above theorems are sharp, in the sense that
when $H=1/2$, we get exactly the same constant as in the inequality
(5.5) of \cite{dgw} with the metric $d_2$. For the $T_1$ inequality,
the sharpness will be discussed in the next section, where we will
apply the above results to study small-time and large-time asymptotics
of the solution of a fractional stochastic differential equation (SDE).

%s3 ###
\section{Small-time and large-time asymptotics of the solution of a
fractional SDE}\label{appli}

The concentration inequalities on the path space of continuous
functions are very well adapted to investigate small- and large-time
asymptotics of processes. The link between the concentration
inequalities and the $L^1$ transpor\-ta\-tion inequality is proved in
\cite{bg}. We recall that a measure $\mu$ on the metric space $(E,d)$
satisfies the property $T_1(C)$ if and only if for any Lipschitzian
function $F \dvtx (E,d)\to\mathbb R$, $F$ is $\mu$-integrable and for all
$\lambda\in\mathbb R,$ we have the Gaussian concentration inequality
\[
\int_E \exp\biggl( \lambda \biggl( F - \int_E F\,\mathrm{d}\mu  \biggr) \biggr)\,\mathrm{d}\mu \le\exp\biggl(C
\| F\|_{\mathrm{Lip}} \frac{ \lambda^2}2 \biggr),
\]
where
\[
\| F\|_{\mathrm{Lip}} = \sup_{x\neq y} \frac{|F(x)-F(y)|}{d(x,y)} .
\]
By Chebyshev's inequality and an
optimization argument, we obtain the following Hoeffding-type
inequality:
\begin{equation}\label{hoeffding}
\mu\biggl( F - \int_E F\,\mathrm{d}\mu > r \biggr) \le\exp\biggl(- \frac{r^2}{2 C  \|
F\|_{\mathrm{Lip}}^2} \biggr)\qquad  \forall r>0 .
\end{equation}
We present Hoeffding-type inequalities for the solution $X$ of
\eqref{eds} on the metric space of continuous functions associated with
the metrics $d_\infty$ and $d_2$.

Let $V\dvtx R^d \to\mathbb R$ be a function such that $\| V\|_{\mathrm{Lip}} \le
\alpha$. We consider $F$ and $F_\infty$ defined\break on $C(0,T;\mathbb R^d)$
by
\begin{eqnarray*}
F(\gamma) &=&\frac{1}T\int_0^T V(\gamma(t))\,\mathrm{d}t, \\
F_\infty(\gamma) & = &\sup_{t\in[0,T]} |\gamma(t)-\gamma(0)|.
\end{eqnarray*}
The function $F$ is $\alpha$-Lipschitzian with respect to $d_\infty$
and $\alpha/\sqrt{T}$-Lipschitzian with respect to the metric $d_2$. As
for $F_\infty$, it is $1$-Lipschitzian with respect to the metric
$d_\infty$. The following properties are consequences of
\eqref{hoeffding}.

\subsection*{Small-time asymptotics}\vspace*{-2pt}
There exists a constant $C$ (depending only on $H$ and $\sigma$) such
that
if we assume (H1) (resp., (H2)), then
the solution of \eqref{eqt1-1} (resp., \eqref{eqdim1}) satisfies, for
all $r>0$ and small $T$,
\begin{equation}\label{hoeffding1}
\mathbb P_x \biggl(\frac{1}T \int_0^T  [ V(X_t) -\mathbb E V(X_t)  ]\,\mathrm{d}t
>r \biggr) \le\exp\biggl(- \frac{r^2}{C \alpha^2 T^{2H}} \biggr) ,
\end{equation}
and using \eqref{hoeffding} with the functional $F_\infty$ yields that
there exists some $C$ such that
\begin{equation}\label{hoeffding2}
\mathbb P_x \Bigl(  \Bigl[\sup_{t\in[0,T]}| X_t-x| -\mathbb E \Bigl(\sup_{t\in[0,T]}|
X_t-x|  \Bigr)  \Bigr] >r \Bigr) \le\exp\biggl(- \frac{r^2}{2 C T^{2H}} \biggr) .
\end{equation}

\vspace*{-2pt}\subsection*{Large-time asymptotics}\vspace*{-2pt}
In the framework of Theorem \ref{t2} (resp., Theorem \ref{t3}), we
assume that (H3) (resp., H4(c)) is satisfied for
$B<0$. The solution of equation \eqref{eqt1-1} (resp., equation
\eqref{eqdim1}) satisfies the following: for any
$r>0$,
\begin{eqnarray}\label{hoeffding3}
\mathbb P_x \biggl(
\frac{1}T\int_0^T  [ V(X_t) -\mathbb E V(X_t)  ] \,\mathrm{d}t >r \biggr)
&\le&\exp\biggl(- \frac{r^2 B^2 T^{2-2H} }{4\alpha^2 H \|\sigma\|
_{0,T,\infty}^2( 1-e^{BT} )}\biggr)\\[-2pt]
&&\hspace*{-9.9pt}\mbox{(resp.,} \nonumber\\[-2pt]
\label{hoeffding3bis}
 &\le&\exp\biggl(- \frac{r^2 B^2 T^{2-2H}}{4\alpha^2 H
\sigma_1^2\sigma_2^2( 1-e^{BT/\sigma_1} )} \biggr)).\vspace*{-1pt}
\end{eqnarray}

\begin{remi*}
\begin{enumerate}[(ii)]
\item[(i)] When $H=1/2$, the inequality \eqref{hoeffding2} gives the
correct order when $T\to0+$ (see \cite{dgw}, Remark 5.12(b)). This
justifies that the constants $C$ in the $T_1(C)$
properties established in our work are of correct order and are sharp
in some sense.
\item[(ii)] The estimates \eqref{hoeffding3} and
\eqref{hoeffding3bis} are well adapted to the study of large-time
asymptotics of the solutions of \eqref{eqt1-1} and \eqref{eqdim1}.
These estimates are sharp, in the sense that when we put $H=1/2$ into
the formula, we obtain the same Hoeffding-type estimate as given in
\cite{dgw} (see Corollary 5.11).\vspace*{-2pt}
\end{enumerate}
\end{remi*}
%

%s4 ###
\section{Deterministic differential equations driven by rough
functions}\vspace*{-2pt}\label{prelim}
This section deals with deterministic differential equations driven by
H\"{o}lder continuous functions. These equations are the ones satisfied
by the trajectories of the solution of equation \eqref{eqt1-1}. Our aim
is to prove an estimate with respect to the metric $d_\infty$ for the
difference of two solutions of deterministic differential equations
driven by two different H\"{o}lder continuous functions. This is
clearly the first step if we want to use a Gaussian tail
criterion.\vadjust{\goodbreak}

Suppose that $f\in C^{\lambda}(a,b)$ and $g\in C^{\mu}(a,b)$ with
$\lambda+ \mu>1$. From \cite{young}, the Riemann--Stieltjes integral
$\int_a^b f\,\mathrm{d}g$ exists. In \cite{za}, the author provides an explicit
expression for the integral $\int_a^bf\,\mathrm{d}g$ in terms of fractional
derivatives. Let $\alpha$ be such that $\lambda>\alpha$ and
$\beta>1-\alpha$. Supposing that the following limit exists and is
finite, we define $g_{b-}(t) = g(t)-\lim_{\eps\downarrow0}g(b-\eps)$.
The Riemann--Stieltjes integral can then be expressed as
%e2 ###
%
\begin{equation}\label{a6}
\int_{a}^{b}f_{t}\,\mathrm{d}g_{t}=(-1)^{\alpha}\int_{a}^{b}( D_{a+}^{\alpha }f)
(t)( D_{b-}^{1-\alpha}g_{b-}) (t)\,\mathrm{d}t,
\end{equation}
where
\[
D_{a+}^{\alpha}f(t)=\frac{1}{\Gamma(1-\alpha)}\biggl( \frac{f(t)}{%
(t-a)^{\alpha}}+\alpha\int_{a}^{t}\frac{f(t)-f(s)}{(t-s)^{\alpha
+1}}\,\mathrm{d}s\biggr)
\]
and
\[
D_{b-}^{\alpha}g_{b-}(t)=\frac{(-1)^{\alpha}}{\Gamma(1-\alpha )}\biggl(
\frac{g(t)-g(b)}{(b-t)^{\alpha}}+\alpha\int_{t}^{b}\frac
{g(t)-g(s)}{%
(s-t)^{\alpha+1}}\,\mathrm{d}s\biggr).
\]
We refer to \cite{skm} for further details on fractional operators. We
first state the following useful lemma concerning the estimation of
integrals like \eqref{a6}. The proof is identical to the one proposed
in \cite{hn} and so we only highlight some constants.
\begin{lemme}\label{esti-int}
For $0<\beta<1$ and $f, g$ in $C^\beta(0,T;\mathbb R^d),$ there exists
a constant $\kappa$ such that for any $0\le a <b\le T$,
\begin{equation}\label{oo}
\biggl| \int_{a}^{b}f_{t}\,\mathrm{d}g_{t} \biggr| \le\frac{\kappa}{\beta-1/2} \
\|g\|_{0,T,\beta}  [ \|f\|_{a,b,\infty}(b-a)^\beta+ \|f\|_{a,b,\beta}
(b-a)^{2\beta}  ] .
\end{equation}
\end{lemme}
\begin{pf}
We choose $\alpha$ such that $1-\beta<\alpha<1/2$ and use (\ref{a6})
to write that for all $0\leq s,t\leq T$,
\[
\biggl| \int_s^t f_r\,\mathrm{d}g_r \biggr| \leq\int_s^t  | D^{\alpha}_{s+}f_r
D^{1-\alpha}_{t-}g_{t-}(r)|\,\mathrm{d}r .
\]
We have
\begin{eqnarray*}
 | D^{1-\alpha}_{t-}g_{t-}(r)  | & \leq&
\frac{\beta}{(\alpha+\beta-1)\Gamma(\alpha)} \|g\|_{0,T,\beta}
|t-r|^{\alpha+\beta-1}\quad  \mbox{and}\\
 | D^{\alpha}_{s+}f_r | &\leq&
\frac{\|f\|_{s,t,\infty}}{\Gamma(1-\alpha)} (r-s)^{-\alpha} +
\frac{\alpha\|f\|_{s,r,\beta}}{(\beta-\alpha)\Gamma(1-\alpha)}
(r-s)^{\beta-\alpha} .
\end{eqnarray*}
It follows that
\begin{eqnarray*}
\biggl| \int_s^t f_r\,\mathrm{d}g_r \biggr| & \leq&
\frac{\beta\|f\|_{s,t,\infty}\|g\|_{0,T,\beta}}{(\alpha+\beta
-1)\Gamma(\alpha)\Gamma(1-\alpha)}
\int_s^t (r-s)^{-\alpha}(t-r)^{\alpha+\beta-1}\,\mathrm{d}r \\ %
&&{} + \frac{\beta\alpha\|f\|_{s,t,\beta}\| g\|_{0,T,\beta}
}{(\beta-\alpha)(\alpha+\beta-1)\Gamma(\alpha)\Gamma(1-\alpha)}
\int_s^t (r-s)^{\beta-\alpha}(t-r)^{\alpha+\beta-1}\,\mathrm{d}r .
\end{eqnarray*}
We use the change of variables $r=(t-s)\xi+s$ and, recalling that the
beta function is defined by $\mathcal B(a,b)= \int_{0}^1
(1-\xi)^{a-1}\xi^{b-1}\,\mathrm{d}\xi= \frac{\Gamma(a)\Gamma(b)}{\Gamma(a+b)}$,
we get
\[
\hspace*{-50.2pt}\biggl| \int_s^t f_r \,\mathrm{d}g_r \biggr| \leq k_{\alpha,\beta} \| g\| _{0,T,\beta} [
\|f\|_{s,t,\infty} (t-s)^\beta+ \|f\|_{s,t,\beta} (t-s)^{2\beta}]
\]
with
\begin{eqnarray*}
\hspace*{30pt}k_{\alpha,\beta} & = &\frac{\beta\mathcal
B(\alpha+\beta,1-\alpha)}{(\alpha+\beta-1)\Gamma(\alpha)\Gamma
(1-\alpha)} +\frac{\alpha\beta
\mathcal B(\alpha+\beta,1+\beta-\alpha)}{(\alpha+\beta
-1)(\beta-\alpha)\Gamma(\alpha)\Gamma(1-\alpha)}\\[-2pt]
& \le&\frac{\kappa}{\beta-1/2} := c_\beta .
\end{eqnarray*}
The fact that $k_{\alpha,\beta} \leq\kappa/(\beta-1/2),$ where
$\kappa$ is a universal constant independent of $\alpha$ and~$\beta$,
is proved in \cite{s}.
\end{pf}

Set $1/2 <\beta<1$ and let $g,\tilde g\in C^\beta(0,T;\R^m)$. We shall
work with two deterministic differential equations on $\rd$:
\begin{eqnarray*}
x_{t}^{i} &=& x_{0}^i + \int_{0}^t b^i (x_s)\,\mathrm{d}s + \sum_{j=1}^m \int_{0}^t
\s^{i,j}(s)\,\mathrm{d}g_{s}^{j},\qquad t\in\oT, \\[-2pt]
\tilde x_{t}^{i} &=& x_{0}^i + \int_{0}^t b^i (\tilde x_s)\,\mathrm{d}s +
\sum_{j=1}^m \int_{0}^t \s^{i,j}(s)\,\mathrm{d}\tilde g_{s}^{j},\qquad  t\in \oT ,
\end{eqnarray*}
$i=1,\ldots,d$, $x_0\in\rd$.

It is proved in \cite{nr}, Theorem 5.1 that if $1-\beta<\alpha< 1/2$,
then each of the above equations has a unique $(1-\alpha)$-H\"{o}lder
continuous solution. The estimates on the solution $(x_t)_{t\in[0,T]}$
obtained in \cite{nr} were improved in \cite{hn}, Theorem 3.3.
Unfortunately, these estimates are unusable in our context.
Nevertheless, since the matrix $\sigma$ does not depend on the
solution, our framework is more simple, and we quickly prove the
estimate we need in the following proposition.
\begin{prop}\label{th-esti}
Let $g$ and $\tilde g$ be H\"{o}lder continuous of order $1/2<\beta<1$.
Under the assumptions~(\textup{H1}), we define $\Delta=
(2L_b)^{-1}\wedge1$. For all $T\le\Delta$, there exists a universal
constant $K$ such that
\[
\|x-\tilde x\|_{0,T,\infty} \le K\| \s\|_\beta\| g-\tilde g\|
_{0,T,\beta} T^\beta .
\]
\end{prop}
\begin{pf} We restrict ourselves to the case $d=m=1$ for
simplicity. We write
\[
x_t-\tilde x_t = \int_0^t  [ b(x_r)- b(\tilde x_r) ]\,\mathrm{d}r +\int_0^t \s(r)\,\mathrm{d}[ g_r -\tilde g_r ] .
\]
Using \eqref{oo},
we may write
\[
|x_t-\tilde x_t | \leq t  L_b  \| x-\tilde x\|_{0,t,\infty} +
c_\beta\| g-\tilde g\|_{0,t,\beta} \|\s\|_{\beta} [ t^\beta+
t^{2\beta}],\vadjust{\goodbreak}
\]
where $c_\beta= \kappa/(\beta-1/2),$ and consequently
\[
\|x-\tilde x\|_{0,t,\infty}  \leq t L_b  \| x-\tilde x\|_{0,t,\infty}
+ c_\beta\| g-\tilde g\|_{0,t,\beta} \|\s\|_{\beta} [ t^\beta+
t^{2\beta}].
\]
Therefore the result is proved when $t\le\Delta$.
\end{pf}

%
%s5 ###
\section{Proofs of the main results}\label{preuves}
%s5.1 ###
\subsection{$T_1(C)$ for paths of SDE's driven by an fBm}
To prove Theorem \ref{t1}, we use a sufficient condition that is
present in the proof of \cite{dgw}, Theorem 2.3. This is recalled in
the following lemma whose proof is entirely contained in the
aforementioned proof.
\begin{lemme}
Let $\mu$ a probability measure on a metric space $(E,d)$. Let $\xi$
and $\xi'$ be two independent random variables valued in $E$ with law
$\mu$ defined on some probability space $(\Omega,\mathcal F,\mathbb
P)$. If
\[
C :=2\sup_{k\ge1} \biggl( \frac{k!  \mathbb E (d(\xi,\xi'))^{2k}}{(2k)!}\biggr)^{1/k}
\]
is finite, then $\mu$ satisfies the transportation inequality $T_1(C)$
on $(E,d)$.
\end{lemme}
We now turn to the proof of Theorem \ref{t1} itself.
\begin{pf*}{Proof of Theorem \ref{t1}}
Let $(B_t^H)_{t\in[0,T]}$ and $(\tilde B_t^H)_{t\in[0,T]}$ be two
independent fractional Brownian motions defined on the filtered
probability space $(\Omega,\mathcal F,(\mathcal F_t)_{t\in[0,T]},
\mathbb P)$. We denote by $(X_t)_{t\in[0,T]}$ and $(\tilde X_t)_{t\in
[0,T]}$ the strong solutions of \eqref{eqt1-1} driven by $B$ and
$\tilde B$, respectively. The $T_1(C)$ property will be implied by the
finiteness of
\[
C = 2\sup_{k\ge1} \biggl( \frac{k!  \mathbb E  ( d^{2k}_{\infty} (
X,\tilde X ) )}{(2k)!} \biggr)^{1/k} .
\]
Let $1/2<\beta< H<1$ and $T\le\Delta$. Proposition \ref{th-esti}
implies that
\[
d^{2k}_{\infty}  ( X,\tilde X ) \le K^{2k} \|\s\|_{\beta }^{2k}
\|B^H-\tilde B^H\|_{0,T,\beta}^{2k}  T^{2k\beta} .
\]
In the following, the constant $K$ is universal, but may vary from line
to line. Taking expectation and using \eqref{esti-mom} from Lemma
\ref{lemme-grr}, we obtain
\begin{eqnarray*}
C & \le& 2\sup_{k\ge1} \biggl( \frac{k!  K^{2k} \|\s\|_{\beta }^{2k}\
T^{2k\beta}  T^{2k(H-\beta)} (2k)! }{k! (2k)!}
\biggr)^{1/k}\\
& \le& K \|\s\|_\beta T^{2H} ,
\end{eqnarray*}
and the result is proved.\vadjust{\goodbreak}~%
\end{pf*}
\begin{pf*}{Proof of Theorem \ref{t1-dim1}}
If we set
\[
F(y) = \int_0^y \frac{\mathrm{d}z}{\sigma(z)} ,
\]
then we can use
the change-of-variables formula \cite{za}, Theorem 4.3.1 to obtain that
$(X_t)_{t\in[0,T]}$ is the unique solution of
\[
X_t  = x +\int_0^t b(X_s)\,\mathrm{d}s +\int_0^t \sigma(X_s)\,\mathrm{d}B_s^H,\qquad 0\le t\le T,
\]
if and only if the process $(Y_t)_{t\in[0,T]}$ defined by $Y_t =
F(X_t)$ is the unique solution of
\begin{equation}\label{eqz}
Y_t = F(x) + \int_0^t \frac{b(F^{-1}(Y_s))}{\sigma(F^{-1}(Y_s))}\,\mathrm{d}s
+B_t^H,\qquad  0\le t\le T .
\end{equation}
Our result will follow from the stability of the transportation
inequalities under a Lipschitzian map (see \cite{dgw}, Lemma 2.1). We
consider the map $\Psi$ from the metric space $(C(0,T),d_\infty)$ into
itself defined by $\Psi(\gamma) = F^{-1}\circ\gamma$. We have, for
$\gamma_1,\gamma_2 \in C(0,T),$
\[
d_\infty \bigl( \Psi(\gamma_1) -\Psi(\gamma_2)  \bigr) \le\|\Psi'\| _\infty
d_\infty(\gamma_1,\gamma_2)
\]
and, clearly, $\Psi' = (F^{-1})' = \sigma$.
Thus, the map $\Psi$ is $\alpha$-Lipschit\-zian with
$\alpha=\sigma_2$. If $\mathbb P_x^X$ (resp., $\mathbb P^Y_{F(x)}$)
denotes the law of the process $X$ (res\-p., $Y$), then
\[
\mathbb P^X_x=\mathbb P^Y_{F(x)} \circ F = \mathbb P^{Y}_{F(x)} \circ
\Psi^{-1} .
\]
We denote by $L_{\tilde b}$ the Lipschitz constant of the function
$\tilde b = b\circ F^{-1} / \sigma\circ F^{-1}$. It is easy to check
that
\[
L_{\tilde b} \le\frac{\sigma_2}{\sigma_1^2}  ( L_b \sigma_2 + L_\sigma
B ).
\]
By Theorem \ref{t1}, $\mathbb P^Y_{F(x)} \in T_1(K T^{2H})$ for
$T\le(2L_{\tilde b})^{-1}\wedge1$, so we have that $\mathbb P^X_x \in
T_1(K  \sigma_2^2  T^{2H})$ for $T\le\tau$ with
\[
\tau= 1 \wedge\frac{\sigma_1^2 }{ 2 \sigma_2 ( L_b \sigma_2 +
L_\sigma B )}.\vspace*{9pt}
\]
\upqed\end{pf*}

%s5.2 ###
\subsection{$T_2(C)$ for paths of SDE's driven by an fBm}\vspace*{6pt}
In this subsection, we prove Theorems \ref{t2} and \ref{t3}. First, we
briefly recall some basic facts about stochastic integration with
respect to fBm. We refer to \cite{n} for a more detailed
treatment.

\subsubsection*{Preliminaries}\vspace*{-2pt}
Let $\mathcal H$ be the Hilbert space defined as the closure of
$\mathcal E$ (the set of step functions on $\oT$ with values in
${\mathbb R}^m$)
with respect to the scalar product\vspace*{-2pt}
\[
\bigl\langle\bigl(\1_{[0,t_1]},\ldots,\1_{[0,t_m]}\bigr),
\bigl(\1_{[0,s_1]},\ldots,\1_{[0,s_m]}\bigr) \bigr\rangle_{\mathcal H} = \sum_{i=1}^m
R_H(t_i,s_i).\vspace*{-2pt}
\]
%
%Then the scalar product between two elements $\ffi$ and $\psi$ of
%$\mathcal E$ is given by
%& \langle\ffi,\psi\rangle_{\mathcal H} = H(2H-1)
%|r-u|^{2H-2} \ffi^i_r\psi^i_udrdu .%
%The space $\mathcal H$ contains ${L}^{\frac{1}H}(0,T;\r^m)$ but
%its element can be distributions. Formula
%
The mapping $(\1_{[0,t_1]},\ldots,\1_{[0,t_m]}) \mapsto\sum_{i=1}^m
B^{H,i}_{t_i}$ is extended to an isometry between $\mathcal H$ and the
Gaussian space $H_1(B^H)$ associated with $B^H$. We denote this
isometry by $\ffi\mapsto B(\ffi)$.
%
%The covariance kernel $R_H$ can be written as
%$$ R_H(t,s)= \int_{0}^{s\wedge t} K_H(s,u)K_H(t,u)du  .$$
Using the kernel $K$ defined in \eqref{defK}, we introduce the operator
$\mathcal K_{H}^{\ast}\dvtx\mathcal H\to
L^2(0,T;{\mathbb R}^m)$:\vspace*{-2pt} %
\begin{equation}\label{khetoile}
( \mathcal{K}_{H}^{\ast}\varphi) (s) =\int
_{s}^{T}\varphi(r)%
\frac{\partial K_{H}}{\partial r}(r,s)\,\mathrm{d}r .\vspace*{-2pt}
\end{equation}
We have $ \mathcal K_{H}^{\ast}(  ( \1_{[0,t_1]},\ldots,\1_{[0,t_m]} ) ) =
( K_H(t_1,\cdot),\ldots,K_H(t_m,\cdot))$ and, for $\ffi,\psi\in\mathcal E,$\vspace*{-2pt}
\[
\langle\ffi,\psi\rangle_{\mathcal H} = \langle \mathcal K_{H}^{\ast}
\ffi,\mathcal K_{H}^{\ast} \psi\rangle_{L^{2}(0,T;{\mathbb
R}^m)}=\mathbb E( B^H(\ffi )B^H(\psi) ).\vspace*{-2pt}
\]
$\mathcal K_{H}^{\ast}$ then provides an isometry between the
Hilbert space $\mathcal H$ and a closed subspace of $L^{2}(0,T;{\mathbb
R}^m)$. %

We have already mentioned the transfer principle (see \eqref{transfer})
when $B^H$ is written as an integral of the underlying Brownian motion
$W$. More precisely, the transfer principle means that
%
%. Additionally, the Brownian motion $W$ can be expressed as
%W_t = B^H ((\mathcal K_H^\ast)^{-1} (\1_{\ot},...,\1_{\ot}
%) ) .
for any $\ffi\in\mathcal H$, $B^H(\ffi) = W( \mathcal K_H^\ast \ffi)$.

We define $\mathcal{K}_{H}\dvtx L^{2}(0,T;\mathbb{R}^{m})\rightarrow
\mathcal{H}%
_{H}:=\mathcal{K}_{H}( L^{2}(0,T;\mathbb{R}^{m}))$, the operator
defined by $\mathcal K_H h = (\mathcal K_H h^1,\ldots,\mathcal K_H h^m)$
with\vspace*{-2pt}
\[
(\mathcal{K}_{H}h^i)(t):=\int_{0}^{t}K_{H}(t,s)h^i(s)\,\mathrm{d}s,\qquad i=1,\ldots,m.\vspace*{-2pt}
\]
We will use of the following property \cite{d}, Lemma 3.2: for $h\in
L^2(0,T;\mathbb R^m)$,\vspace*{-2pt}
\begin{equation}\label{kh-hol}
| (\mathcal{K}_{H}h)(t) -(\mathcal{K}_{H}h)(s)| \le c |t-s|^H \|
h\|_{L^2(0,T;\mathbb R^m)} .\vspace*{-2pt}
\end{equation}
Using Fubini's theorem and the fact that $\frac{\partial K_H} {\partial
u}(u,s) = c_H( \frac us )^{H-1/2} (u-s) ^{H-3/2}$, we obtain
that if $f\in C^{\lambda}(0,T)$ with $\lambda+H>1$ and $\rho\in
L^2(0,T),$ then it holds that\vspace*{-2pt}
%e3 ###
%
\begin{equation}\label{f4}
\int_{0}^{T}f(r)\,\mathrm{d}(\mathcal{K}_{H}\rho)_{r}=\int_{0}^{T}f(r)\biggl(
\int_{0}^{r}\frac{\partial K_{H}}{\partial r}(r,t)\rho(t)\,\mathrm{d}t\biggr)\,\mathrm{d}r .\vspace*{-2pt}
\end{equation}
The integral on the left-hand side of \eqref{f4} is a
Riemann--Stieltjes integral for H\"{o}lder functions (see Section
\ref{prelim}).

Finally, if $\ffi,\psi\in L^2 (0,T;\mathbb R^m)$, then the scalar
product on $\mathcal H$ has the integral form\vspace*{-2pt}
\[
\langle\ffi{,} \psi\rangle= H(2H-1) \int_0^T\int_0^T |s-t|^{2H-2}
\langle\ffi(s){,} \psi(t)\rangle_{\mathbb R^m}\,\mathrm{d}s\,\mathrm{d}t\vspace*{-2pt}
\]
and, consequently, for $\ffi\in L^2 (0,T;\mathbb R^m),$ we have\vspace*{-2pt}
\begin{equation}\label{inegH}
\| \varphi\|_{\mathcal H}^2 \le 2H T^{2H-1} \| \varphi
\|_{L^2(0,T;\mathbb R^m)}^2 .\vspace*{-2pt}\vadjust{\goodbreak}
\end{equation}

\begin{pf*}{Proof of Theorem \ref{t2}}
We recall that a classical $m$-dimensional Brownian motion $(W_t)_{t\in
[0,T]}$ is defined on $(\Omega, \mathcal F, \mathbb P)$ and
$B^H=(B_t^H)_{t\in[0,T]}$ is an $m$-dimensional fBm defined on
$(\Omega,\mathcal F,\mathbb P)$ transferred from $W$. Let $\mathbb Q$
be a probability measure on $C(0,T;\mathbb R^d)$ such that $\mathbb Q
\ll\mathbb P_x$. We can assume that $\mathbf{H}(\mathbb
Q|\mathbb P_x)<\infty,$ otherwise there is nothing to prove.

The first part of the proof follows the arguments of \cite{dgw}. The
idea is to express the finiteness of the entropy by means of the energy
of the drift arising from the Girsanov transform of a well-chosen
probability measure. This method also appears in \cite{fu} and was well
known for a long time. The relationship between the finite entropy
condition and the finite energy condition on the Girsanov drift
appeared in \cite{follmer1,follmer2} for the first time (to the best of
our knowledge) in the particular case of Brownian motion with
drift.\looseness=1

We consider
\[
\tilde{\mathbb Q} = \frac{\mathrm{d}\mathbb Q}{\mathrm{d}\mathbb P_x}(X) \mathbb P.
\]
Clearly, $\tilde{\mathbb Q}$ is a probability measure on
$(\Omega,\mathcal F)$ and
\begin{eqnarray*}
\mathbf{H}(\tilde{\mathbb Q} | \mathbb P) & =& \int_{\Omega} \ln\biggl(
\frac{\mathrm{d}\tilde{\mathbb Q}}{\mathrm{d}\mathbb P} \biggr)\,\mathrm{d} \tilde {\mathbb Q}\\
& =&
\int_{\Omega} \ln\biggl( \frac{\mathrm{d}\mathbb Q}{\mathrm{d}\mathbb P_x}(X)
\biggr) \frac{\mathrm{d}\mathbb Q}{\mathrm{d}\mathbb P_x}(X)\,\mathrm{d}\mathbb P \\
& = &\int_{C(0,T;\mathbb R^d)} \ln\biggl( \frac{\mathrm{d}\mathbb Q}{\mathrm{d}\mathbb P_x} \biggr)
\frac{\mathrm{d}\mathbb Q}{\mathrm{d}\mathbb P_x}\,\mathrm{d}\mathbb P_x =
\mathbf{H}(\mathbb Q|\mathbb P_x).
\end{eqnarray*}
Following \cite{dgw}, there exists a predictable process
$\rho=(\rho^1(t),\ldots,\rho^m(t))_{0\le t\le T}$ such that
\[
\mathbf{H}(\mathbb Q|\mathbb P_x)=
\mathbf{H}(\tilde{\mathbb Q} | \mathbb P) = \frac{1}2  \mathbb
E_{\tilde{\mathbb Q}} \int_0^T |\rho(t)|^2\,\mathrm{d}t
\]
and, by Girsanov's theorem, the process $(\tilde B_t)_{t\in[0,T]}$
defined by
\[
\tilde B_t = W_t - \int_0^t \rho(s)\,\mathrm{d}s
\]
is a Brownian motion under $\tilde{\mathbb Q}$ and is associated
(thanks to the transfer principle) with the $\tilde{\mathbb
Q}$-fractional Brownian motion $(\tilde{B}^H)_{t\in[0,T]}$ defined by
\[
\tilde B^H_t  = \int_0^t K_H(t,s)\,\mathrm{d}\tilde B_s = \int_0^t K_H(t,s)\,\mathrm{d}W_s
-  ( \mathcal K_H \rho )(t)= B^H_t -  ( \mathcal K_H \rho )(t).
\]
Consequently, under $\tilde{\mathbb Q}$, $X$ verifies
\begin{equation}\label{eqtilde}
\cases{
\mathrm{d}X_t = b(X_t)\,\mathrm{d}t + \sigma(t)\,\mathrm{d}\tilde B^H_t + \sigma(t)\,\mathrm{d}(\mathcal K_H
\rho)(t), \cr
X_0 = x. }
\end{equation}
We now consider the solution $Y$ (under $\tilde{\mathbb Q}$) of the
following equation:
\begin{equation}\label{eqy}
\cases{
\mathrm{d}Y_t = b(Y_t)\,\mathrm{d}t + \sigma(t)\,\mathrm{d}\tilde B^H_t\vspace*{2pt} \cr
Y_0 =x.}
\end{equation}
Under $\tilde{\mathbb Q}$, the law of the process $(Y_t)_{t\in[0,T]}$
is exactly $\mathbb P_x$. Then, $(X,Y)$ under $\tilde{\mathbb Q}$ is a~coupling of $(\mathbb Q,\mathbb P_x)$ and it follows that
\begin{eqnarray*}
[W_2^{d_2}(\mathbb Q,\mathbb P_x)]^2  &\le&  \mathbb E_{\tilde{\mathbb Q}}
( | d_2 ( X,Y ) |^{2}  )=\mathbb E_{\tilde{\mathbb Q}}\biggl( \int_0^T
|X_t-Y_t|^2\,\mathrm{d}t \biggr),\\ [2pt]
[W_2^{d_{\infty}}(\mathbb Q,\mathbb P_x)]^2  &\le&  \mathbb E_{\tilde{\mathbb Q}} ( | d_\infty( X,Y ) |^{2}  )=\mathbb
E_{\tilde{\mathbb Q}}\Bigl(\sup_{0\le t\le T} |X_t-Y_t|^2 \Bigr).
\end{eqnarray*}
We now estimate the distance on $C(0,T;\mathbb R^m)$ between $X$ and
$Y$ with respect to the distances $d_2$ and $d_\infty$. We note that
equations \eqref{eqtilde} and \eqref{eqy} can be considered as pathwise
integral equations driven by $\beta$-H\"{o}lder functions with
$\beta<H$. Indeed, the H\"{o}lder regularity is straightforward for the
driving function $\tilde B$ since it is a fractional Brownian motion
under $\tilde{\mathbb Q}$ (and so it has almost surely
$\beta$-H\"{o}lder trajectories for any $\beta<H$). Moreover, since
$\int_0^T |\rho(s)|^2\,\mathrm{d}s < +\infty$ almost surely, $\mathcal K_H \rho
\in C^H(0,T)$ almost surely by~\eqref{kh-hol}.\looseness=1

We write
\[
X_t-Y_t = \int_0^t  \bigl( b(X_s)-b(Y_s) \bigr)\,\mathrm{d}s +\int_0^t \sigma (s)\,\mathrm{d}
(\mathcal K_H\rho ) (s) .
\]
We use the change of variables formula for a $\beta$-H\"{o}lder
continuous function (see \cite{za}, Theorem~4.3.1) and the stability
assumption (H2) to obtain
\begin{eqnarray}\label{ui}
|X_t-Y_t|^2 & = &2 \sum_{i=1}^d\sum_{j=1}^m \int_0^t  ( X^i_s-Y^i_s
 )\sigma^{i,j}(s)\,\mathrm{d}(\mathcal K_H\rho)^j(s) \nonumber\\[2pt]
&&{}+2 \int_0^t  \langle X_s-Y_s , b(X_s)-b(Y_s)
 \rangle_{\mathbb R^d}\,\mathrm{d}s \\[2pt]
& \le &2 \sum_{i=1}^d\sum_{j=1}^m \int_0^t  ( X^i_s-Y^i_s  )
\sigma^{i,j}(s)\,\mathrm{d}(\mathcal K_H\rho^j)(s) +2B \int_0^t |X_s-Y_s|^2\,\mathrm{d}s.\nonumber
\end{eqnarray}
Since $X-Y\in C^{\beta}(0,T;\mathbb R^d)$ and $\rho\in
L^2(0,T;\mathbb
R^m)$, we use \eqref{khetoile} and \eqref{f4} to obtain
\begin{eqnarray*}
&&\int_0^t  ( X^i_s-Y^i_s  ) \sigma
^{i,j}(s)\,\mathrm{d}(\mathcal K_H\rho^j)(s) \\[2pt]
&&\quad = \int_0^t  ( X^i_s-Y^i_s  ) \sigma^{i,j}(s)\biggl( \int
_0^s\frac{\partial K_H}{\partial s}(s,r) \rho^j(r)\,\mathrm{d}r\biggr)\,\mathrm{d}s\\[2pt]
&&\quad = \int_0^t \biggl(\int_r^t  ( X^i_s-Y^i_s  ) \sigma ^{i,j}(s)
\frac{\partial K_H}{\partial s}(s,r)\,\mathrm{d}s \biggr) \rho^j(r)\,\mathrm{d}r\\[2pt] %
&&\quad = \int_0^t \mathcal K^{\ast}_H  \bigl( (X^i-Y^i )\sigma^{i,j}\1_{[0,t]} \bigr)
(r) \rho^j(r)\,\mathrm{d}r  .
\end{eqnarray*}
We denote by $\sigma^\ast$ the transpose matrix of $\sigma$ and we use
the inequality \eqref{inegH} to obtain
\begin{eqnarray*}
&& 2 \sum_{i=1}^d\sum_{j=1}^m \int_0^t  ( X^i_s-Y^i_s
)\sigma^{i,j}(s)\,\mathrm{d}(\mathcal K_H\rho)^j(s) \\[-1pt]
&&\quad = 2 \int_0^t \bigl\langle\mathcal K^{\ast}_H   \bigl(
\sigma^\ast(X-Y )\1_{[0,t]} \bigr) (r) {,} \rho(r) \bigr\rangle
_{\mathbb R^m}\,\mathrm{d}r \\[-1pt]
&&\quad  \le2 \bigl\|\mathcal K^{\ast}_H  \bigl( \sigma ^\ast(X-Y
)\1_{[0,t]} \bigr) \bigr\|_{L^2(0,T)}   \| \rho\|
_{L^2(0,t)} \\[-1pt]
&&\quad  \le2 \| \sigma^\ast(X-Y )\1_{[0,t]} \|_{\mathcal
H}   \| \rho\|_{L^2(0,t)} \\[-1pt]
&&\quad  \le2 (2H)^{1/2} T^{H-1/2} \bigl\| \sigma^\ast(X-Y )\1
_{[0,t]} \bigr\|_{L^2(0,T)}   \| \rho\|_{L^2(0,t)} \\[-1pt]
&&\quad  \le2 (2H)^{1/2} T^{H-1/2}\| \sigma\|_{0,T,\infty} \|
X-Y \|_{L^2(0,t)}   \| \rho\|_{L^2(0,t)}  .
\end{eqnarray*}
We report this estimate in \eqref{ui}, and using the inequality
$4\epsilon ab \le4\epsilon^2 a^2 +b^2$ with $\epsilon= ( H T^{2H-1}
\|\sigma\|_{0,T,\infty}^2/(2|B|))^{1/2}$, we obtain
\begin{eqnarray*}
|X_t-Y_t|^2 & \le& 2 (2H)^{1/2} T^{H-1/2} \|\sigma\|_{0,T,\infty}
\| X-Y \|_{L^{2}(0,t)}  \| \rho\|_{L^2(0,t)} \\[-1pt] %
& &{}+2B \int_0^t | X_s-Y_s|^2 \,\mathrm{d}s \\[-1pt]
& \le&(2/|B|)HT^{2H-1} \|\sigma\|_{0,T,\infty}^2 \int_0^t |\rho
(s)|^2\,\mathrm{d}s \\[-1pt] %
& &{} +(2B+|B|) \int_0^t | X_s-Y_s|^2\,\mathrm{d}s .
\end{eqnarray*}
Gronwall's lemma implies that for any $t>0,$
\[
|X_t-Y_t|^2 \le(2/|B|)HT^{2H-1} \|\sigma\|_{0,T,\infty}^2\int_{0}^t
\mathrm{e}^{(2B+|B|)\times(t-s)}|\rho(s)|^2\,\mathrm{d}s.
\]
Hence, we may write that
\[
d_\infty^2 (X{,}Y) \le(2H/|B|)T^{2H-1}
\|\sigma\|_{0,T,\infty}^2 \bigl(1\vee \mathrm{e}^{(2B+|B|)\times T}\bigr)\int_0^T
|\rho(s)|^2\,\mathrm{d}s
\]
and
\[
[W_2^{d_\infty}(\mathbb Q,\mathbb P_x)]^2 \le2 C_{T,H}
\mathbf{H}(\mathbb Q | \mathbb P_x)
\]
with $C_{T,H}=2H T^{2H-1}(1\vee \mathrm{e}^{(2B+|B|)\times
T})\|\sigma\|_{0,T,\infty}^2 /|B|$.\vadjust{\goodbreak}

Analogously for the metric $d_2$, we have
\begin{eqnarray*}
[W_2^{d_2}(\mathbb Q,\mathbb P_x)]^2 & \le&\mathbb
E_{\tilde{\mathbb Q}}\int_0^T | X_t-Y_t|^2\,\mathrm{d}t \\ %
& \le &(2/|B|)HT^{2H-1}\|\sigma\|_{0,T,\infty}^2 \\ %
& &{} \times\mathbb E_{\tilde{\mathbb Q}} \int_{0}^T
|\rho(s)|^2 \biggl( \int_{s}^T \mathrm{e}^{(2B+|B|)\times(t-s)}\,\mathrm{d}t \biggr)\,\mathrm{d}s. %
\end{eqnarray*}
Since
\[
\int_{s}^T \mathrm{e}^{(2B+|B|)\times(t-s)}\,\mathrm{d}t \le\cases{
\displaystyle\frac{\mathrm{e}^{3BT}-1}{3B}, &\quad if $B> 0$,\vspace*{2pt} \cr
\displaystyle-\frac{1-\mathrm{e}^{BT}}{B}, &\quad if $B<0$,
}
\]
we define
\[
c_{B,T} := \cases{
\displaystyle\frac{\mathrm{e}^{3BT}-1}{3}, &\quad if $B> 0$,\vspace*{2pt} \cr
1-\mathrm{e}^{BT}, &\quad if $B<0$
}
\]
and it follows that
\begin{eqnarray*}
[W_2^{d_2}(\mathbb Q,\mathbb P_x)]^2 & \le & 4
(H/B^2)T^{2H-1}\|\sigma\|_{0,T,\infty}^2 c_{B,T} \biggl( \frac{1}2 \mathbb
E_{\tilde{\mathbb Q}} \int_{0}^T |\rho(s)|^2\,\mathrm{d}s\biggr)\\
& \le & 2C_{T,H}  \mathbf{H}(\mathbb Q | \mathbb P_x)
\end{eqnarray*}
with $C_{T,H} = (2/B^2)HT^{2H-1}\|\sigma\|_{0,T,\infty}^2 c_{B,T}$.
\end{pf*}
\begin{pf*}{Proof of Theorem \ref{t3}}
We use the same change-of-variables as in the proof of Theorem~\ref{t1-dim1} and
consider the map $\Psi$ from the metric space
$(C(0,T),d_2)$ into itself defined by $\Psi(\gamma) = F^{-1}\circ
\gamma$. We have, for $\gamma_1,\gamma_2 \in C(0,T),$
\begin{eqnarray*}
d_2  \bigl( \Psi(\gamma_1) -\Psi(\gamma_2)  \bigr) & =& \biggl( \int
_0^T |\Psi(\gamma_1(s)) - \Psi(\gamma_2(s))|^2 \,\mathrm{d}s \biggr)^{1/2} \\
& \le & \|\Psi'\|_\infty d_2(\gamma_1,\gamma_2),
\end{eqnarray*}
thus the map $\Psi$ is $\sigma_2$-Lipschitzian. If $\mathbb P_x^X$
(resp., $\mathbb P^Y_{F(x)}$) denotes the law of the process $X$
(resp.,~$Y$), then
\[
\mathbb P^X_x=\mathbb P^Y_{F(x)} \circ F = \mathbb
P^{Y}_{F(x)} \circ\Psi^{-1}.
\]
Since $\mathbb P^Y_{F(x)} \in T_2(C)$,
we have that $\mathbb P^X_x \in T_2(\sigma_2^2 C)$. It remains to
prove that the stability assumption (H3) is true for the
function $\tilde b = b\circ F^{-1} / \sigma\circ F^{-1}$. Writing
$\tilde b' = ( b'\circ F^{-1} \sigma\circ F^{-1}-b\circ F^{-1}
\sigma'\circ F^{-1})/ \sigma\circ F^{-1}$, it easy to see that under
the\vadjust{\goodbreak} assumptions (H4), we have\looseness=1
\begin{eqnarray*}
\biggl(x-y , \frac{(b\circ F^{-1})(x)}{(\sigma\circ F^{-1})(x)}-
\frac{(b\circ F^{-1})(y)}{(\sigma\circ F^{-1})(y)}\biggr)  \le
\frac{B}{\sigma_1} |x-y|^2 .
\end{eqnarray*}\looseness=0
We can then apply Theorem \ref{t2} to equation \eqref{eqz} and thus the
result (b) is proved. A~similar reasoning is true for the
metric $d_\infty$.
\end{pf*}

%s6 ###
\section{A remark on the link between the exponential moment and~$T_1(C)$}\label{surprise}

It has been proven in \cite{bg,bv,dgw} that $\mu\in T_1(C)$ if and
only if we have, for some $\delta>0,$ the Gaussian tail
\[
\int_E \int_E \mathrm{e}^{\delta d^2(x,y)} \mu(\mathrm{d}x)\mu(\mathrm{d}y) <+\infty .
\]
The link between the constant $C$ and the exponential moment is
described in the following remark.
\begin{remi*}\label{mom}
Let $\mu$ a probability measure on a metric space $(E,d)$. Assume that
there exists some $\delta>0$ such that the following Gaussian tail
holds:
\[
C(\delta) :=\int_E \int_E \mathrm{e}^{\delta d^2(x,y)} \mu(\mathrm{d}x)\mu(\mathrm{d}y)
<+\infty.
\]
Then, $\mu$ satisfies the transportation inequality $T_1(C)$ on
$(E,d)$. In \cite{dgw}, the authors have linked $C$ with the above
exponential moment in the following way:
\begin{equation}\label{link-c}
C \le\frac{2}\delta\sup_{k\ge1} \biggl( \frac{(k!)^2}{(2k)!} \int_E \int_E
\mathrm{e}^{\delta d^2(x,y)} \mu(\mathrm{d}x)\mu(\mathrm{d}y) \biggr)^{1/k} .
\end{equation}
By an optimization argument, the supremum in the formula \eqref{link-c}
is achieved for $k=1$ and consequently $C \le C(\delta)/\delta$.
\end{remi*}
In \cite{bv} (see also \cite{goz}, page 69), the authors have proven
that the constant $C$ is in fact controlled by a better constant, but
it is not tractable to study short-time and long-time asymptotic
behavior.

In our context, if we use the above remark and the exponential estimate
\eqref{esti-exp} from Lemma \ref{lemme-grr}, then we can easily prove
that the law $\mathbb P_x$ of the solution of equation \eqref{eqt1-1}
satisfies the property $T_1(C)$ with $C=K \|\sigma\|_\beta T^{2H-\eps}$
for small time $T$ and a small $\eps>0$. Nevertheless, the power of $T$
is not the correct order when one applies this result to small-time
asymptotics.

We believe that it remains an interesting open problem to give a simple
link between the exponential moment and the constant $C$ in $T_1(C)$.
\begin{pf*}{Proof of the estimate $C\le C(\delta)/\delta$}
We use an optimization argument involving the gamma function $\Gamma$.
We denote, for $x\ge1,$
\[
\Phi(x)  = \exp\biggl( \frac{1}x \ln\biggl( C(\delta)\frac {\Gamma^2
(x+1)}{\Gamma(2x+1)} \biggr) \biggr) .
\]
We remark that the right-hand side of \eqref{link-c} is equal to
$(2/\delta) \Phi(k)$. Our result will then be a~consequence of
$\sup_{x\ge1} \Phi(x) = \Phi(1) = C(\delta)/2$. We denote by $\Psi$
the function $(\ln\Gamma)'=\Gamma'/\Gamma$ (usually called the
\emph{digamma function}). We write $\Phi'(x) = h(x) \Phi(x)/x^2$,
where the function~$h$ is defined for $x\ge1$ by
\[
h(x)  = - \ln\biggl( C(\delta)\frac{\Gamma^2 (x+1)}{\Gamma(2x+1)} \biggr) +2x
 \bigl(\Psi(x+1)-\Psi(2x+1) \bigr).
\]
Obviously, $\Phi'$ and $h$ have the same sign. Since $\Psi'(x) =
\sum_{k=0}^\infty\frac{1}{(x+k)^2}$ (see \cite{aar}, page 13), we
deduce that
\begin{eqnarray*}
&&\Psi'(x+1)-2\Psi'(2x+1)\\
& &\quad= \sum_{k=0}^\infty\frac{1}{(x+1+k)^2} -\frac
{1}{2(x+(k+1)/2)^2} \\
& &\quad=\frac{1}2 \sum_{k=0}^\infty\frac{1}{(x+1+k)^2}
+\frac{1}2\sum_{k=0}^\infty\frac{1}{(x+1+k)^2} -\frac
{1}{(x+(k+1)/2)^2} \\
& &\quad = \frac{1}2\Biggl\{\sum_{k=0}^\infty\frac {1}{(x+1+k)^2}
+\sum_{j=0}^\infty-\frac{1}{(x+(2j+1)/2)^2}\Biggr\}
\\
& &\quad = \frac{1}2\Biggl\{\sum_{k=0}^\infty\frac
{1}{(x+1+k)^2} +\sum_{k=0}^\infty-\frac{1}{(x+k+1/2)^2}\Biggr\} \\
& &\quad= \frac{1}2 \{ \Psi'(x+1)-\Psi'(x+1/2)\}.
\end{eqnarray*}
Since $\Psi''(x) = -2\sum_{k=0}^\infty(x+k)^{-3}$, $\Psi'$ is a
decreasing function and then
\[
\Psi'(x+1)-2\Psi'(2x+1)\le0 .
\]
This
yields $h'(x) = 2x (\Psi'(x+1)-2\Psi'(2x+1))\le0$. So, for any $x\ge
1$,
\[
h(x)\le h(1)= -\ln\bigl( C(\delta)\Gamma^2(2)/\Gamma(3) \bigr) +2 \bigl(
\Psi(2)-\Psi(3) \bigr) .
\]
In \cite{aar}, the following identity is stated
for $n\ge1$:
\[
\Psi(x+n) = \sum_{k=0}^{n-1} \frac{1}{x+k} +\Psi(x),
\]
so $\Psi(2)-\Psi(3) =-1/2$. Finally, $h(1)= -\ln(C(\delta)/2)-1\le0$
because $C(\delta)\ge1$. Thus, $h(x)\le0$ for any $x\ge1,$ and
$\Phi$ is decreasing. Its maximum is achieved for $x=1$.
\end{pf*}

%%%%%%%%%%%%%%%%%%%%%%%%%%%%%%%%%%%%%%%%%%%%%%%
%%%%%%%%%%%%%%%%%%%%%%%%%%%%%%%%%%%%%%%%%%%%%%%%
%%%%%%%%%%%%%%%%%%%%%%%%%%%%%%%%%%%%%%%%%%%%%%%%
%%%%%%%%%%%%%%%%%%%%%%%%%%%%%%%%%%%%%%%%%%%%%%%%
%%%%%%%%%%%%%%%%%%%%%%%%%%%%%%%%%%%%%%%%%%%%%%%%%%%%%%%%%%%%%%%%%%%%%%%
%%%%%%%%%%%%%%%%%%%%%%%%%%%%%%%%%%%%%%%%%%%%%%%%%%%%%%%%%%%%%%%%%%%%%%%
%%%%%%%%%%%%%%%%%%%%%%%%%%%%%%%%%%%%%%%%%%%%%%%%%%%%%%%%%%%%%%%%%%%%%%%
%%%%%%%%%%%%%%%%%%%%%%%%%%%%%%%%%%%%%%%%%%%%%%%%%%%%%%%%%%%%%%%%%%%%%%%%%%%%%%%%%%%%%
\begin{appendix}\label{appendix}
%s7 ###
\section*{Appendix: Fernique-type lemma}
\renewcommand{\theequation}{\arabic{equation}}
\setcounter{equation}{21}
\begin{lemme}\label{lemme-grr}
Let $T>0$, $1/2< \beta<H <1$. Then, for any $\alpha<
1/(128(2T)^{2(H-\beta)}),$
\begin{equation}\label{esti-exp}
\mathbb E  [ \exp( \alpha\|B^H\|_{0,T,\beta}^2 )  ]
\leq\bigl(1-128\alpha(2T)^{2(H-\beta)}\bigr)^{-1/2} .
\end{equation}
Moreover, we have the following moment estimate for any $k\ge1$:
\begin{equation}\label{esti-mom}
\mathbb E (\| B^H\|_{0,T,\beta}^{2k} ) \le32^k (2T)^{2k(H-\beta)}
\frac{(2k)!}{k!} .
\end{equation}
\end{lemme}
\begin{pf}
First, we prove that
%e4 ###
%
\begin{equation}\label{modul}
|B^{H,i}_t - B^{H,i}_s | \leq\xi_{\beta}  |t-s|^{\beta},\qquad i=1,\ldots,m,
\end{equation}
where $\xi_{\beta}$ is a positive random variable such that
\begin{equation}\label{modul1}
E( \xi_{\beta}^{2p} )\leq32^p (2T)^{2p(H-\beta)} \frac{(2p)!}{p!} .
\end{equation}
Although the proofs of \eqref{modul} and \eqref{modul1} are classical,
we include them for the convenience of the reader. With $\psi(u)=
u^{2/(H-\beta)}$ and $p(u)=u^{H}$ in Lemma 1.1 of \cite{grr}, the
Garsia--Rodemich--Rumsey inequality reads as follows:
\[
|B^{H,i}_t-B^{H,i}_s| \leq8 \int_0^{|t-s|}\biggl( \frac{4 \Delta}{u^2}
\biggr)^{(H-\beta)/2} H u^{H-1}\,\mathrm{d}u,
\]
where the random variable $\Delta$ is
\[
\Delta= \int_0^T\int_0^T \frac{|B^{H,i}_t-B^{H,i}_s|^{2/(H-\beta
)}}{|t-s|^{2H/(H-\beta)}}\,\mathrm{d}t\,\mathrm{d}s .
\]
We have
\begin{eqnarray*}
|B^{H,i}_t-B^{H,i}_s| & \leq& 8   (4\Delta)^{(H-\beta)/2}  \int
_0^{|t-s|}H  u^{\beta-1}\,\mathrm{d}u \leq8   (4\Delta)^{(H-\beta)/2}\
\frac{H}{\beta}  |t-s|^{\beta} \\ %
& \leq& 8   (4\Delta)^{(H-\beta)/2} |t-s|^{\beta}  .
\end{eqnarray*}
We let $\xi_{\beta} = 8   (4\Delta)^{(H-\beta)/2}$ and for $p\ge
1/(H-\beta), $ we have
\begin{eqnarray*}
\E\xi_{\beta}^{2p}
& \leq& 8^{2p} 4^{p(H-\beta)} \mathbb E \biggl( \int_0^T\int_0^T
\frac{|B^{H,i}_t-B^{H,i}_s|^{2/(H-\beta)}}{|t-s|^{2H/(H-\beta)}}
\,\mathrm{d}t\,\mathrm{d}s \biggr)^{p(H-\beta)} \\
& \leq& 8^{2p} (2T)^{2p(H-\beta)} \int_0^T\int_0^T \frac{\mathbb E
|B^{H,i}_t-B^{H,i}_s|^{2p}}{|t-s|^{2pH}} \frac{\mathrm{d}t\,\mathrm{d}s}{T^2}\\
& \leq& 8^{2p} (2T)^{2p(H-\beta)}\frac{(2p)!}{2^p p!} \leq32^p
(2T)^{2p(H-\beta)}\frac{(2p)!}{p!}  .
\end{eqnarray*}
Thus, \eqref{modul} and \eqref{modul1} are proved. What remains to be
shown can be tediously deduced from \cite{fer}, Theorem 1.3.2. We can
also make the following direct computations. Using \eqref{modul} and
\eqref{modul1}, we have
\begin{eqnarray*}
\mathbb E ( \exp ( \alpha\| B^H \|_{\beta}^2  )  ) & \leq& \mathbb E (
\exp ( \alpha\xi_{\beta}^2  ) ) \leq\mathbb E\Biggl(
\sum_{p=0}^\infty\frac{\alpha^p \xi_{\beta}^{2p}}{p!} \Biggr)\\
& \leq & \sum_{p=0}^{\infty} (32\alpha)^p
(2T)^{2p(H-\beta)}\frac{(2p)!}{(p!)^2} \\
& \le & \bigl( 1-128\alpha(2T)^{2(H-\beta)}\bigr)^{-1/2} ,
\end{eqnarray*}
where we have used the identity $\sum_{p=0}^\infty
a^p\frac{(2p)!}{(p!)^2} = (1-4a)^{-1/2}$ for $a<1/4$. Thus, the lemma
is proved.
\end{pf}
\end{appendix}

\section*{Acknowledgements} I am grateful to Francis Hirsch for his
valuable comments as well as to an anonymous referee for his careful
reading.

\printhistory

\end{document}